\documentclass{notices}
\usepackage{amsfonts,amssymb,amsmath,amscd,amsthm}

\usepackage{xcolor}
\usepackage{bm}
\usepackage{comment}
\usepackage{mathtools}
\usepackage{empheq} 
\mathtoolsset{showonlyrefs}

\newcommand{\yx}[1]{\textcolor{blue}{\textsf{[YX: #1]}}}

\newcommand{\db}[1]{\textcolor{green}{\textsf{[DB: #1]}}}

\newcommand{\ud}{\,\mathrm{d}}
\newcommand{\PP}{\mathbb{P}}
\newcommand{\RR}{\mathbb{R}}
\newcommand{\mc}[1]{\mathcal{#1}}
\newcommand{\bX}{\bm{X}}
\newcommand{\bl}{\bm{\ell}}
\newcommand{\bh}{\bm{h}}
\newcommand{\bx}{{\bm{x}}}

\newcommand{\nn}{\mathcal{NN}}
\newcommand{\half}{\frac{1}{2}}
\newcommand{\transpose}{^{\operatorname{T}}}
\newcommand{\bW}{\bm{W}}
\newcommand{\balpha}{\bm{\alpha}}

\newcommand{\ltwonorm}[1]{\left\lVert#1\right\rVert_2}

\DeclareMathOperator*{\argmin}{arg\,min}

\newtheorem{defi}{Definition}

\title{Pandemic Control, Game Theory and Machine Learning}

\author{
  Yao Xuan
  \affil{Yao Xuan contributed to the project when he was a Ph.D. student at the University of California, Santa Barbara. His email address is yxscience@gmail.com.}
  \and
  Robert Balkin
  \affil{Robert Balkin is a current graduate student supervised by Ruimeng Hu and Hector D. Ceniceros at the University of California, Santa Barbara. His email address is rbalkin@ucsb.edu.}
  \and Jiequn Han
   \affil{Jiequn Han is a research fellow at the Flatiron Institute. His email address is jiequnhan@gmail.com.}
  \and Ruimeng Hu
   \affil{Ruimeng Hu is an assistant professor at the University of California, Santa Barbara. Her email address is rhu@ucsb.edu.}
   \and Hector D. Ceniceros
   \affil{Hector D. Ceniceros is a full professor at the University of California, Santa Barbara. His email address is ceniceros@ucsb.edu.}
}

\begin{document}

\maketitle

\section*{COVID-19 and Control Policies}
The coronavirus disease 2019 (COVID-19) pandemic has brought an enormous impact on our lives. Based on data from World Health Organization, as of May 2022, there have been more than 520 million confirmed cases of infection and more than 6 million deaths globally; In the United States, there have been more than 83 million confirmed cases of infection and more than one million cases of death. Needless to say,  the economic impact has also been catastrophic, resulting in unprecedented unemployment and the bankruptcy of many restaurants, recreation centers, shopping malls, etc. 

Control policies play a crucial role in the alleviation of the COVID-19 pandemic. For example, lockdown and work-from-home policies and mask requirements on public transport and public areas have been proved to be effective in stopping the spreading of COVID-19. On the other hand, governors also have to be aware of the economic activity loss due to these pandemic control policies. Therefore, a thorough understanding of the evolution of COVID-19 and the corresponding decision-making provoked by such a virus will be beneficial for future events and in other interconnected systems around the world.

\subsection*{Epidemiology}

Epidemiology is the science of analyzing the distribution and determinants of health-related states and events in specified populations. It is also the application of this study to the control of health problems. Infectious diseases are one of this kind, including the ongoing novel coronavirus (COVID-19). 

Since March 2020, when the World Health Organization declared the COVID-19 outbreak a global pandemic, epidemiologists have made tremendous efforts to understand how COVID-19 infections emerge and spread and how they may be prevented and controlled. Many epidemiological methods involve mathematical tools, e.g., using causal inference to identify causative agents and factors for its propagation, and molecular methods to simulate disease transmission dynamics.


The first epidemic model concerning epidemic spreading dates back to 1760 by Daniel Bernoulli \cite{bernoulli1760essai}. Since then, many papers have been dedicated to this field and, later on, to epidemic control. Among control strategies, the quarantine, firstly introduced in 1377 in Dubrovnik on Croatia’s Dalmatian Coast \cite{grmek1997beginnings}, has shown as a powerful component of the public health response to emerging and reemerging infectious diseases. However, quarantine and other measures for controlling epidemic diseases have always been controversial due to the potentially raised political, ethical, and socioeconomic issues. Such complication naturally calls for the inclusion of decision-making in epidemic control, as it helps to answer how to take \emph{optimal} actions to balance public interest and individual rights. But not until recent years have there been some research studies in this direction. Moreover, when multiple authorities are involved in the decision-making process, it is challenging to analyze how to collectively or competitively make decisions due to the difficulty of solving this high-dimensional problem.
 


In this article, we focus on the decision-making development for the intervention of COVID-19, aiming to provide mathematical models and efficient numerical methods, and justifications for related policies that have been implemented in the past and explain how the authorities' decisions affect their neighboring regions from a game theory viewpoint.

\subsection*{Mathematical models} 
In a classic, compartmental epidemiological model, each individual in a geographical region is assigned a label, e.g.,
\textbf{S}usceptible, \textbf{E}xposed, \textbf{I}nfectious, \textbf{R}emoved, \textbf{V}accinated. 
Different labels represent different status -- \textbf{S}: those who are not
yet infected; \textbf{E}: who have been infected but are not yet infectious themselves; \textbf{I}: who have been
infected and are capable of spreading the disease to those in the susceptible category, \textbf{R}: who
have been infected and then removed from the disease due to recovery or death, and \textbf{V}: who have been vaccinated and are immune to the infection. 
As COVID-19 progressed, it was learned that spread from asymptomatic cases was an important driving force. More refined models may further split \textbf{I}
into mild-symptomatic/asymptomatic individuals who are in-home for recovery and serious-symptomatic ones that need hospitalization. We point to \cite{MR4121530} which considers a similar problem in the optimal control setting, which includes asymptomatic individuals and the effect of impulses. 



\begin{figure*}[h]
\begin{center}
  \includegraphics[width=\textwidth]{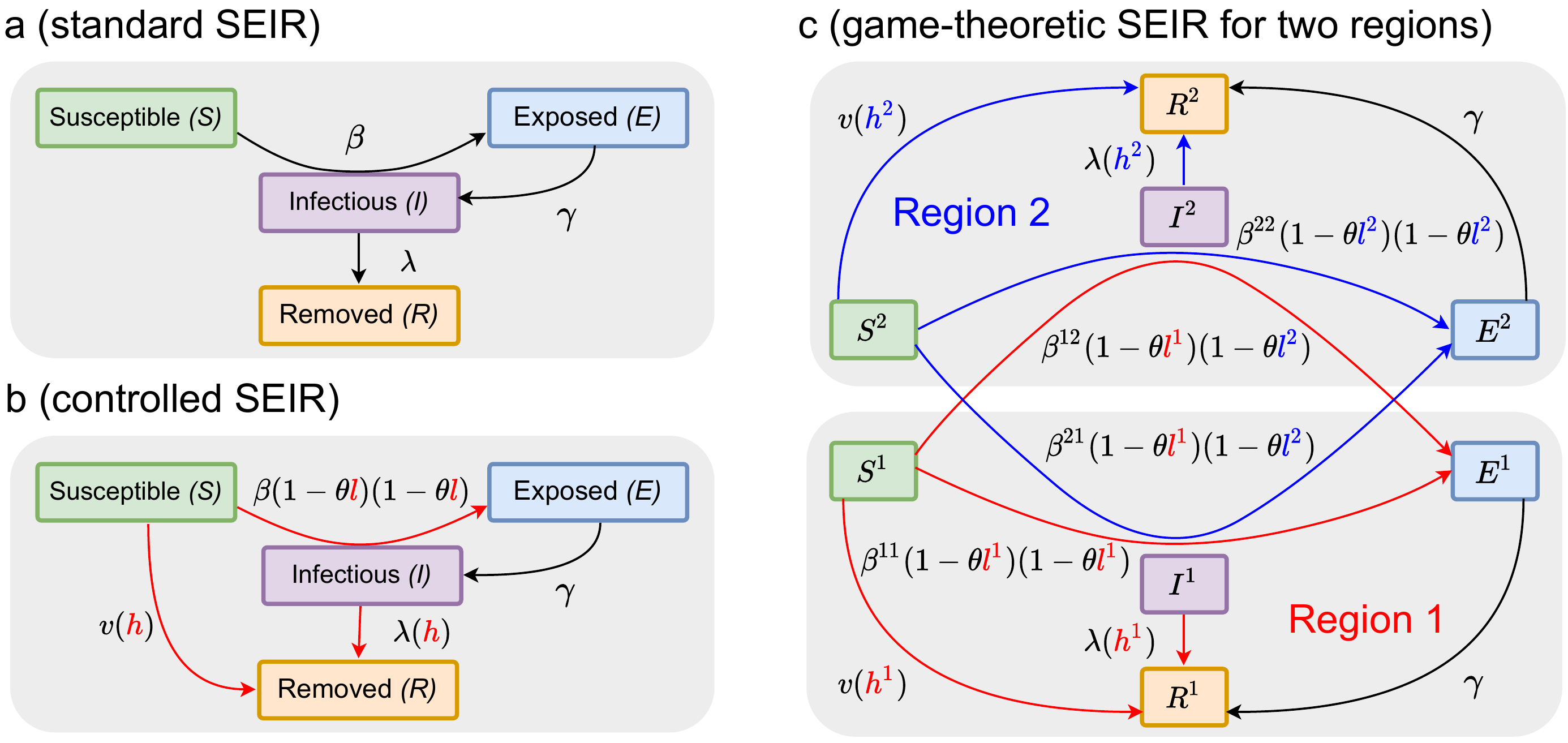}
  \end{center}
  \caption{(a) A simple SEIR model:  susceptible individuals become exposed after close contact with infected ones; those exposed become infectious after a latency period; and those infected become removed afterward due to recovery or death; (b) Controlled SEIR model: the planner chooses the level of nonpharmaceutical policies (lockdown or work from home) $\ell$ and pharmaceutical policies (effort of vaccination development or distribution) $h$ affecting the transitions such that only $(1-\theta \ell(t))$ of the original susceptible and infectious individuals can contact each other, and affecting the recovery rate $\lambda(h)$ from infectious individuals to removed ones, here $\theta$ is used describe the effectiveness of policy $\ell$; (c) An illustration of the game-theoretic SEIR model for two regions. }\label{fig:SEIR}
\end{figure*}

Individuals transit between these compartments, and the labels' order in a model indicates the flow patterns between the compartments. For instance, in a simple SEIR model \cite{MR908379} (see also Figure~\ref{fig:SEIR}a), a susceptible becomes exposed after close contact with infected ones; exposed individuals become infectious after a latency period; and infected ones become removed afterward due to recovery or death. Let $S(t)$, $E(t)$, $I(t)$ and $R(t)$ be the proportion of population of each compartment at time $t$, the following differential equations provide the mathematical model: 
\begin{equation}\label{def:SIR}
\begin{aligned}
 & \dot S(t) =  -\beta S(t)I(t), \\
    &  \dot E(t) = \beta S(t)I(t) - \gamma E(t), \\
    & \dot I(t) = \gamma E(t) - \lambda I(t), \; \dot R(t) = \lambda I(t),
\end{aligned}
\end{equation}
where $\beta$ is  the average number of contacts per person per time, $\gamma$ describes the latent period when the person has been infected but not yet infectious, and $\lambda$ represents the recovery rate measuring the proportion of people recovered or dead from infected population.


Many infections, such as measles and chickenpox, confer long-term, if not lifelong, immunity, while others, such as influenza, do not. As evidenced by numerous epidemiological and clinical studies analyzing possible factors for COVID reinfections, COVID-19 falls precisely into the second category \cite{nordstrom2022risk}. Mathematically, this can be taken into account by adding a transition $I \to S$.


Though deterministic models such as \eqref{def:SIR} have received more attention in the literature, mainly due to their tractability, stochastic models have some advantages. The epidemic-spreading progress is by nature stochastic. Moreover, introducing stochasticity to the system 
could account for numerical and empirical uncertainties, and also provide probabilistic predictions, i.e., a range of possible scenarios associated with their likelihoods. This is crucial for understanding the uncertainties in the estimates.  


One class of stochastic epidemic models uses continuous-time Markov chains, where the state process takes discrete values but evolves in continuous time and is Markovian. In a simple Stochastic SIS (susceptible-infectious-susceptible) model \cite{MR1025386} with a population of $N$ individuals, let $X_t$ be the number of infected individuals at time $t$, $\beta$ the rate of infected individuals infecting those susceptible, and $\lambda$ the rate that an infected individual recovers and becomes susceptible again. The transition probabilities among states $n$, $n+1$, $n-1$ are
\begin{align}
    & \PP(X_{t+ \Delta t} = n+1 \vert X_t = n) \approx \frac{\beta}{N}n(N-n) \Delta t, \\
    & \PP(X_{t+ \Delta t} = n-1 \vert X_t = n) \approx \lambda n \Delta t, \\
    & \PP(X_{t+ \Delta t} = n \vert X_t = n ) \approx 1 - \left(\frac{\beta}{N}n(N-n)  + \lambda n\right)\Delta t.
\end{align}
Another way to construct a stochastic model is by introducing white noise $W_t$ in \eqref{def:SIR} \cites{tornatore2005stability,MR2428373}, which we shall mainly consider in this paper and describe in details in the later section. 




\subsection*{Control of disease spread}


After modeling how diseases are transmitted through a population, epidemiologists then design corresponding control measures and recommend health-related policies to the region planner. 

\smallskip
In general, there are two types of interventions: pharmaceutical interventions (PIs), such as getting vaccinated and taking medicines, and nonpharmaceutical interventions (NPIs), such as requiring mandatory social distancing, quarantining infected individuals, and deploying protective resources. For the ongoing COVID-19, intervention policies that have been implemented include, but are not limited to, issuing lockdown or work-from-home policies, developing vaccines, and later expanding equitable vaccine distribution, 
providing telehealth programs, deploying protective resources and distributing free testing kits, educating the public on how the virus transmits, and focusing on surface disinfection.

\smallskip
Mathematically, this can be formulated as a control problem: the planner chooses the level of each policy affecting the transitions in \eqref{def:SIR} such that the region's overall cost is minimized. Generally, NPIs help mitigate the spread by lowering the infection rate $\beta$, e.g., a lockdown or work-from-home policy $\ell(t)$ implemented at time $t$ modifies the transition to 
\begin{equation}
    \dot S(t) = - \beta (1-\theta\ell(t))S(t) (1 - \theta\ell(t))I(t),
\end{equation}
meaning that
only $(1-\theta\ell(t))$ of the original susceptible and infectious individuals can contact each other where $\theta$ describes the effectiveness of $\ell$ \cite{alvarez2020simple} (see Figure~\ref{fig:SEIR}b). PIs such as taking preventive medicines, if available, will also lower the infection rate $\beta$, while using antidotes will increase the recovery rate $\lambda$. The modeling of vaccinations is more complex. Depending on the target disease, it may reduce $\beta$ (less chance to be infected) or increase $\lambda$ (faster recovery). It may even create a new compartment ``Vaccinated'' in which individuals can not be infected and which is an absorbing state if life-long immunity is gained. 

A region planner, taking into account the interventions' effects on the dynamics \eqref{def:SIR}, decides on policy by weighing different costs. These costs may include the economic loss due to decrease in productivity during a lockdown, the economic value of life due to death of infected individuals and other social-welfare costs due to the aforementioned measurements.






\section*{Game-theoretic SEIR Model}

Game theory studies the strategic interactions among rational players and has applications in all fields of social science, computer science, financial mathematics, and epidemiology. A game is non-cooperative if players cannot form alliances or if all agreements need to be self-enforcing. Nash equilibrium is the most common kind of self-enforcing agreement \cite{MR43432}, in which a collective strategy emerges from all players in the game to which no one has an incentive to deviate unilaterally.

Nowadays, as the world is more interconnected than ever before, one region's epidemic policy will inevitably influence the neighboring regions. For instance, in the US, decisions made by the governor of New York will affect the situation in New Jersey, as so many people travel daily between the two states. Imagine that both state governors make decisions representing their own benefits, take into account others' rational decisions, and may even compete for the scarce resources (e.g., frontline workers and personal protective equipment). These are precisely the features of a non-cooperative game. Computing the Nash equilibrium from such a game will provide valuable, qualitative guidance and insights for policymakers on the impact of specific policies. 



We now introduce a multi-region stochastic SEIR model \cite{xuan2020optimal} to capture the game features in epidemic control. We give an illustration for two regions in Figure~\ref{fig:SEIR}c. Each region's population is divided into four compartments: \textbf{S}usceptible, \textbf{E}xposed, \textbf{I}nfectious, and \textbf{R}emoved. Denote by $S^n_t, E^n_t, I^n_t, R^n_t$ the \emph{proportion} of the population in the four compartments of the region $n$ at time $t$. They satisfy the following stochastic differential equations (SDEs), which have included interventions (PIs and NPIs), stochastic factors, and game features,  
\begin{align}
    \ud S^n_t &= -\sum_{k = 1}^N \beta^{nk} S^n_t  I^k_t  (1-\theta \ell^n_t)(1-\theta \ell^k_t) \ud t \\
 &\quad - v(h^n_t)S^n_t \ud t - \sigma_{s_n} S^n_t \ud W^{s_n}_t, \label{def:St}\\
\ud E^n_t &= \sum_{k = 1}^N \beta^{nk} S^n_t  I^k_t  (1-\theta \ell^n_t)(1-\theta \ell^k_t) \ud t \label{def:Et} \\ 
&\quad - \gamma E^n_t \ud t + \sigma_{s_n} S^n_t \ud W^{s_n}_t - \sigma_{e_n} E^n_t \ud W^{e_n}_t,\\
\ud I^n_t  &= (\gamma E^n_t - \lambda(h^n_t)I^n_t) \ud t + \sigma_{e_n} E^n_t \ud W^{e_n}_t,  \label{def:It} \\
\ud R^n_t  &= \lambda(h^n_t)I^n_t \ud t + v(h^n_t) S^n_t \ud t, \label{def:Rt}
\end{align}
where $n \in \mc{N} := \{1, 2, \ldots, N\}$ is the collection of $N$ regions, $W_t$ with different superscripts indicate white noise for a compartment in a specific region,  $\bm\ell_t \equiv (\ell^1_t, \ldots, \ell^N_t)$ and $\bm h_t \equiv (h^1_t, \ldots, h^N_t)$ are NPIs and PIs chosen by the region planners at time $t$.  The planner of region $n$ minimizes its region's cost within a period $[0,T]$:
\begin{equation}\label{def:J}
\begin{aligned}
J^n(\bm \ell, \bm h) := \mathbb{E} \bigg[ \int_0^T e^{-rt} P^n\big[ (S^n_t + E^n_t + I^n_t)  \ell^n_t w \\
+ a(\kappa I^n_t \chi + p I^n_t c)\big]  + e^{-rt}\eta (h^n_t)^2\ud t \bigg].
\end{aligned}
\end{equation}

We explain the model \eqref{def:St}--\eqref{def:J} in details:

\medskip
\noindent\textbf{S.} In \eqref{def:St}, $\beta^{nk}$ denotes the average number of contacts of infected people in region $k$ with susceptible ones in region $n$ per time unit. 
	Although some regions may not be geographically connected, the transmission between the two is still possible due to air travel, but is still less intensive than the transmission within the region, i.e., $\beta^{nk}>0$ and $\beta^{nn} \gg \beta^{nk}$ for all $k \neq n$. The decision for NPIs of region $n$'s planner is given by $\ell^n_t \in [0, 1]$. In particular, it represents the fraction of the population being under NPIs (such as social distancing) at time $t$. We assume that those under interventions cannot be infected. However, the policy may only be partially effective as essential activities (food production and distribution, health, and basic services) have to continue. We use $\theta \in [0,1]$ to measure this effectiveness. The transition rate under policy $\bm \ell$ thus become $\beta^{nk} S^n_t  I^k_t  (1-\theta \ell^n_t)(1-\theta \ell^k_t)$. The case $\theta = 1$ means the policy is fully effective. One can also view $\theta$ as the level of public compliance.
	
	The planner of region $n$ also makes the decision $h^n_t \in [0,1]$. This represents the effort, at time $t$, that the planner puts into PIs. We refer to this term, $h^n_t$, as the \emph{health policy}. It will influence the vaccination availability $v(\cdot)$ and the recovery rate $\lambda(\cdot)$ of this model. $v(h^n_t)$ denotes the vaccination availability of region $n$ at time $t$. In this model, we assume that once vaccinated, the susceptible individuals $v(h^n_t)S^n_t$ become immune to the disease, and join the removed category $R^n_t$. This assumption is not very consistent with COVID-19 but reasonable for a short-term decision-making problem. We model it as an increasing function of $h^n_t$, and if the vaccine has not yet been developed, we can define $v(x) = 0$ for $x \leq \overline h$.

\medskip
\noindent\textbf{E.}
In \eqref{def:Et}, $\gamma$ describes the latent period when the person is infected but is not yet infectious. It is the inverse of the average latent time and we assume $\gamma$ to be identical across all regions. The transition between $E^n$ and $I^n$ is proportional to the fraction of exposed individuals, i.e., $\gamma E^n_t$. 
	
\medskip
\noindent\textbf{I and R.}
In \eqref{def:It} and \eqref{def:Rt}, $\lambda(\cdot)$ represents the  recovery rate. For the infected individuals, a fraction $\lambda(h^n)I^n$ (including both death and recovery from the infection) joins the removed category $R^n$ per time unit.
	The rate is determined by the average duration of infection $D$. We model the duration and the recovery rate related to the health policy $h^n_t$ decided by its planner.
The more effort put into the region (e.g., expanding hospital capacity and creating more drive-thru testing sites), the more clinical resources the region will have and the more resources will be accessible by patients, which could accelerate the recovery and slow down death. The death rate, denoted by $\kappa(\cdot)$, is crucial for computing the cost of the region $n$. 
	
\medskip
\noindent\textbf{Cost.}
In \eqref{def:J}, each region planner faces four types of cost. One is the economic activity loss due to the lockdown policy, where $w$ is the productivity rate per individual, and $P^n$ is the population of the region $n$. The second one is due to the death of infected individuals. Here, $\kappa$ is the death rate which we assume for simplicity to be constant, and $\chi$ denotes the economic cost of each death.  The hyperparameter $a$ describes how planners weigh deaths and infections as compared to other costs. 
	The third one is the in-patient cost, where $p$ is the hospitalization rate, and $c$ is the cost per in-patient per day. The last term $\eta (h^n_t)^2$ quantifies the grants for health policies. We choose a quadratic form so that the function is concave in $h^n_t$. This is to account for the law of diminishing marginal utility: the marginal utility from each additional unit declines as investment increases. All costs are discounted by an exponential function $e^{-rt}$, where $r$ is the risk-free interest rate, to take into account the time preference.
	Note that region $n$'s cost depends on all regions' policies $(\bm \ell, \bm h)$, as $\{I^k, k \neq n\}$ appearing in the dynamics of $S^n$. Thus we write it as $J^n(\bm \ell, \bm h)$. 

\medskip
\noindent The above model \eqref{def:St}--\eqref{def:Rt} is by no doubt a prototype, and one can generalize it by considering reinfections (adding transmission from $R^n$ to $S^n$), asymptomatic population (adding asymptomatic compartment $A^n$), different control policy for $S^n$ and $I^n$ (using $\ell_S$ and $\ell_I$ in \eqref{def:St}--\eqref{def:Et}), different fatality rates for young and elder population (introducing $\kappa_Y$ and $\kappa_E$ in \eqref{def:J}).


\subsection*{Nash equilibria and the HJB system}
As explained above, the interaction between region planners can be viewed as a non-cooperative game, when Nash equilibrium is the notion of optimality.

\begin{defi}
A Nash equilibrium (NE) is a tuple $(\bm \ell^\ast, \bm h^\ast) = (\ell^{1, \ast}, h^{1, \ast}, \ldots, \ell^{N, \ast}, h^{N, \ast} ) \in \mathbb{A}^{N}$ such that $\forall n \in \mc{N}$ and $(\ell^n, h^n) \in \mathbb{A}$, 
\begin{equation}
J^n(\bm \ell^\ast, \bm h^\ast) \leq J^n((\bm \ell^{-n, \ast}, \ell^n), (\bm h^{-n, \ast}, h^n)),
\end{equation}
where $\bm \ell^{-n,\ast}$ represents strategies of players other than the $n$-th one:
\begin{equation}
    \bm\ell^{-n, \ast} := [\ell^{1, \ast}, \ldots, \ell^{n-1, \ast}, \ell^{n+1, \ast}, \ldots, \ell^{N,\ast}] \in \mathbb{A}^{N-1}.
\end{equation}
Here $\mathbb{A}$ denotes the set of admissible strategies for each player and  $\mathbb{A}^N$ is the produce of $N$ copies of $\mathbb{A}$. For simplicity, we have assumed that all players take actions in the same space.
\end{defi}



Under proper conditions, the NE is obtained by solving $N$-coupled Hamilton-Jacobi-Bellman (HJB) equations via dynamic programming \cite[Section~2.1.4]{MR3752669}. To simplify the notation, we concatenate the states into a vector form $\bX_t \equiv [\bm S_t,\bm E_t, \bm I_t]\transpose \equiv [S_t^1, \cdots, S_t^N, E_t^1, \cdots, E_t^N, I_t^1, \cdots, I_t^N]\transpose \in \RR^{3N}$, and denote its dynamics by 
\begin{equation}\label{def:Xt}
    \ud \bX_t = b(t, \bX_t, \bm \ell(t, \bX_t), \bm h(t, \bX_t)) \ud t  + \Sigma(\bX_t)\ud \bm W_t.
\end{equation}
For the sake of simplicity, we omit the actual definition of $b$, $f^n$ and $\Sigma$ and refer \cite{xuan2020optimal} for further details. Let $V^n(t, \bx)$ be the minimized cost defined in \eqref{def:J} if the system starts at $\bm X_t = x$. Then, $V^n$, $n = 1, \ldots, N$ solves
\begin{multline}\label{def:HJB}
    \partial_t V^n + \inf_{(\ell^n, h^n) \in [0,1]^2} H^n(t,\bm x,(\bl, \bh)(t, \bx), \nabla_{\bx} V^n)  \\+ \half \text{Tr}(\Sigma(\bx)\transpose \text{Hess}_{\bx} V^n \Sigma(\bx)) = 0,
\end{multline}
with $V^n(T,\bx) = 0$, where $H^n$ is the usual Hamiltonian defined by 
\begin{equation}\label{def:H}
    H^n(t,\bm x, \bl, \bh, \bm p) =  b(t,\bm x, \bl, \bh) \cdot \bm p + f^n(t, \bm x, \ell^n, h^n).
\end{equation}



\section*{Enhanced Deep Fictitious Play}

Solving for the NE of the game is equivalent to solving the $N$-coupled HJB equations of dimension $(3N+1)$ defined in Equation~\eqref{def:HJB}. Due to the high dimensionality, this is a formidable numerical challenge. We overcome this through a deep learning methodology we call \textit{Enhanced Deep Fictitious Play}, being broadly motivated by the method of fictitious play introduced by Brown \cite{Br:49,MR0056265}. 

\paragraph{Deep Learning.} Deep learning leverages a class of computational models composed of multiple processing layers to learn representations of data with multiple levels of abstraction \cite{lecun2015deep}. Deep neural networks are effective tools for approximating unknown functions  in high-dimensional space. In recent years, we have witnessed noticeable success in a marriage of deep learning and computational mathematics to solve high-dimensional differential equations. Specifically, deep neural networks show strong capability in solving stochastic control and games \cite{han2016deep, MR3847747,huLauriere2022}. Below, we take a simple example to illustrate how a deep neural network is determined for function approximation. 


Suppose we would like to approximate a map $y=f(x)$ by a neural network $\nn(x,\bm w)$ in which one seeks to obtain appropriate parameters of the network, $\bm w$, through a process called training. This consists of minimizing a loss function that measures the discrepancies between the approximation and true values over the so-called training set $\{x_i\}_{i=1}^N$. Such a loss function has the general form
\begin{equation}\label{loss_dl}
    L(\bm w)=\frac{1}{N} \sum_{i=1}^N L_i(f(x_i),\nn(x_i, \bm w))+\lambda  \mc{R}(\bm w),
\end{equation}
where $\mc{R}(\bm w)$ is a regularization term on the parameters. The first term $L_i(f(x_i),\nn(x_i, \bm w))$ ensures that the predictions of $\nn(x_i, \bm w)$ match approximately the true value $f(x_i)$ on the training set $\{x_i\}_{i=1}^N$.  Here, $L_i$ could be a direct distance like the $L^p$ norm or error terms derived from some complex simulations associated with $f(x_i)$ and $\nn(x_i, \bm w)$. The hyperparameter $\lambda$ characterizes the relative importance between the two terms in $L(\bm w)$.
To find an optimal set of parameters $\bm w^\ast$, one solves the  problem of minimizing $L(\bm w)$ by the stochastic gradient descent (SGD) method \cite{MR3797719}. Regarding the architecture of $\nn(x, \bm w)$, there is a wide variety of choices depending on the problem, for example fully connected neural networks, convolutional neural networks \cite{lecun1995convolutional}, recurrent neural networks \cite{cho2014learning}, 
and transformers \cite{vaswani2017attention}. 
In this work, we choose fully connected neural networks to approximate the solution and constructed the loss function by simulating the backward differential equations corresponding to the HJB equations.

\paragraph{Enhanced Deep Fictitious Play.}

\begin{figure}[h]
\begin{center}
  \includegraphics[width=\columnwidth]{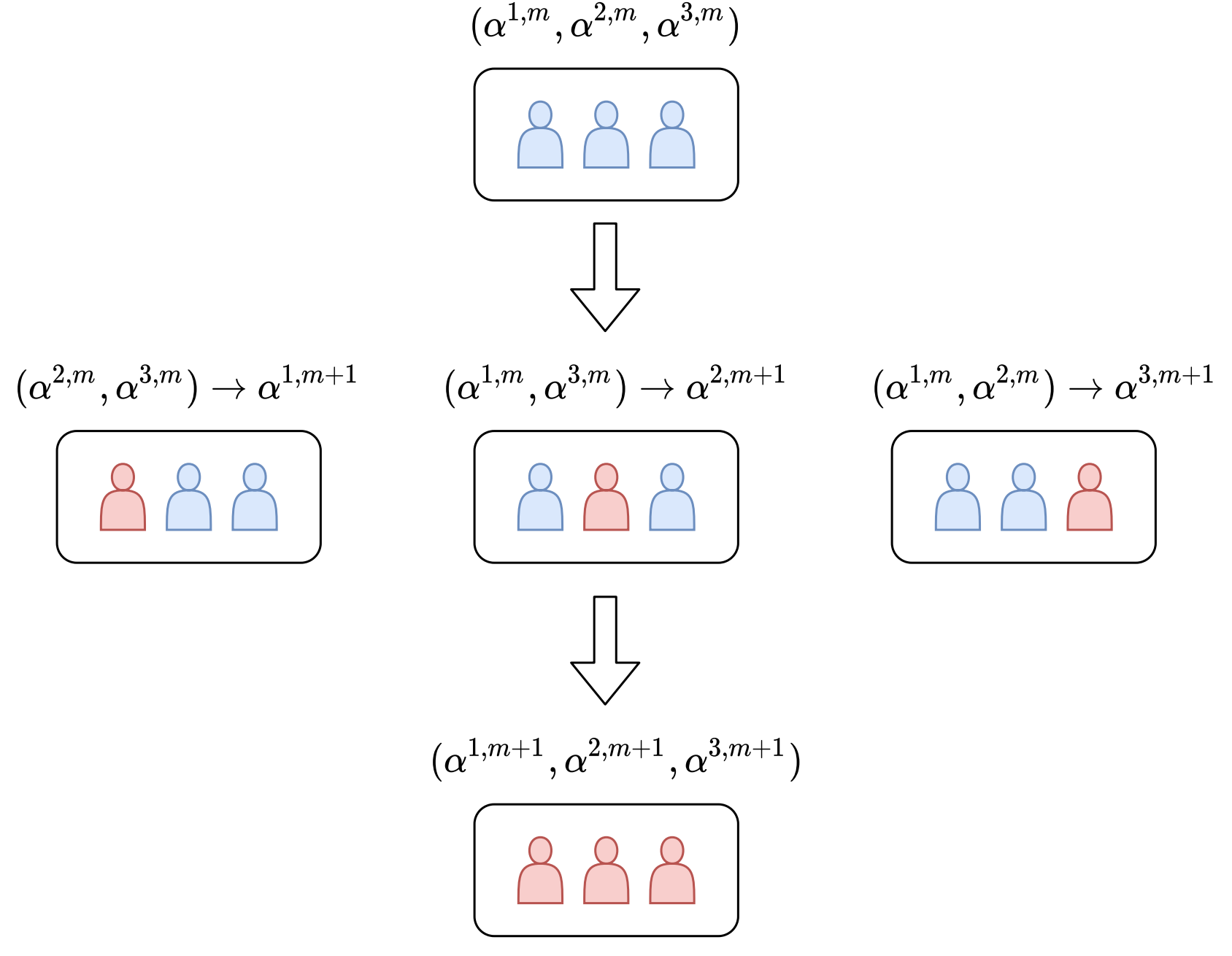}
    \end{center}
  \caption{Schematic plot of fictitious play: each player derives optimal policies at stage $m+1$ assuming other players take optimal strategies at stage $m$. }\label{fp_iterative}
\end{figure}

Note that the HJB system \eqref{def:HJB} is difficult to solve due to the high dimensionality of the $N$-coupled equations. What if we could decouple the system to $N$ separate equations, each of which is easier to solve?   This is the central idea of  \textit{fictitious play},  where we update our approximations to the optimal policies of each player iteratively stage by stage. In each stage, instead of updating the approximations of all the players together by solving the giant system, we do it separately and parallelly. Each player solves for her own optimal policy assuming that the other players are taking their approximated optimal strategies from the last stage. Let us denote the optimal policy and corresponding value function of the single player $n$ in stage $m$ as $\alpha^{n, m}$ and $V^{n, m}$, respectively, and
the collection of these two quantities for all the players as $\bm{\alpha}^{m}=(\alpha^{1,m},...,,\alpha^{N,m})$  and $\bm{V}^{m}=(V^{1,m},...,V^{N,m})$. Finally, let us denote the optimal policies and corresponding value functions for all the players except for player $n$ as $\bm{\alpha}^{-n,m}=(\alpha^{1,m},...,\alpha^{n-1,m},\alpha^{n+1,m},...,,\alpha^{N,m})$ and $\bm{V}^{-n,m}=(V^{1,m},...,V^{n-1,m},V^{n+1,m},...,,V^{N,m})$, where $\alpha^{n,m}$ is a concatenation of lockdown policies and vaccination policies, $i.e., (\ell^{n,m},h^{n,m})$. At stage $m+1$, we can solve for the optimal policy and value function of player $n$  given other players are taken the known policies $\bm{\alpha}^{-n,m}$ and the corresponding value $\bm{V}^{-n,m}$. The logic of fictitious play is shown in Figure~\ref{fp_iterative}, where players iteratively decide optimal policies in stage $m+1$, based on other players' optimal policies in stage $m$. This is slightly different than the usual simultaneous fictitious play, where the belief is described by the time average of past play and the distinction is further discussed in \cite{HaHu:19}. 

The \textit{Enhanced Deep Fictitious Play (DFP)} algorithm we have designed, built from the Deep Fictitious Play (DFP) algorithm \cite{HaHu:19}, reduces time cost from $\mc{O}(M^2)$ to $\mc{O}(M)$ and memory cost from $\mc{O}(M)$ to $\mc{O}(1)$, with $M$ as the total number of fictitious play iterations.

\begin{figure}[h]
  \includegraphics[width=\columnwidth]{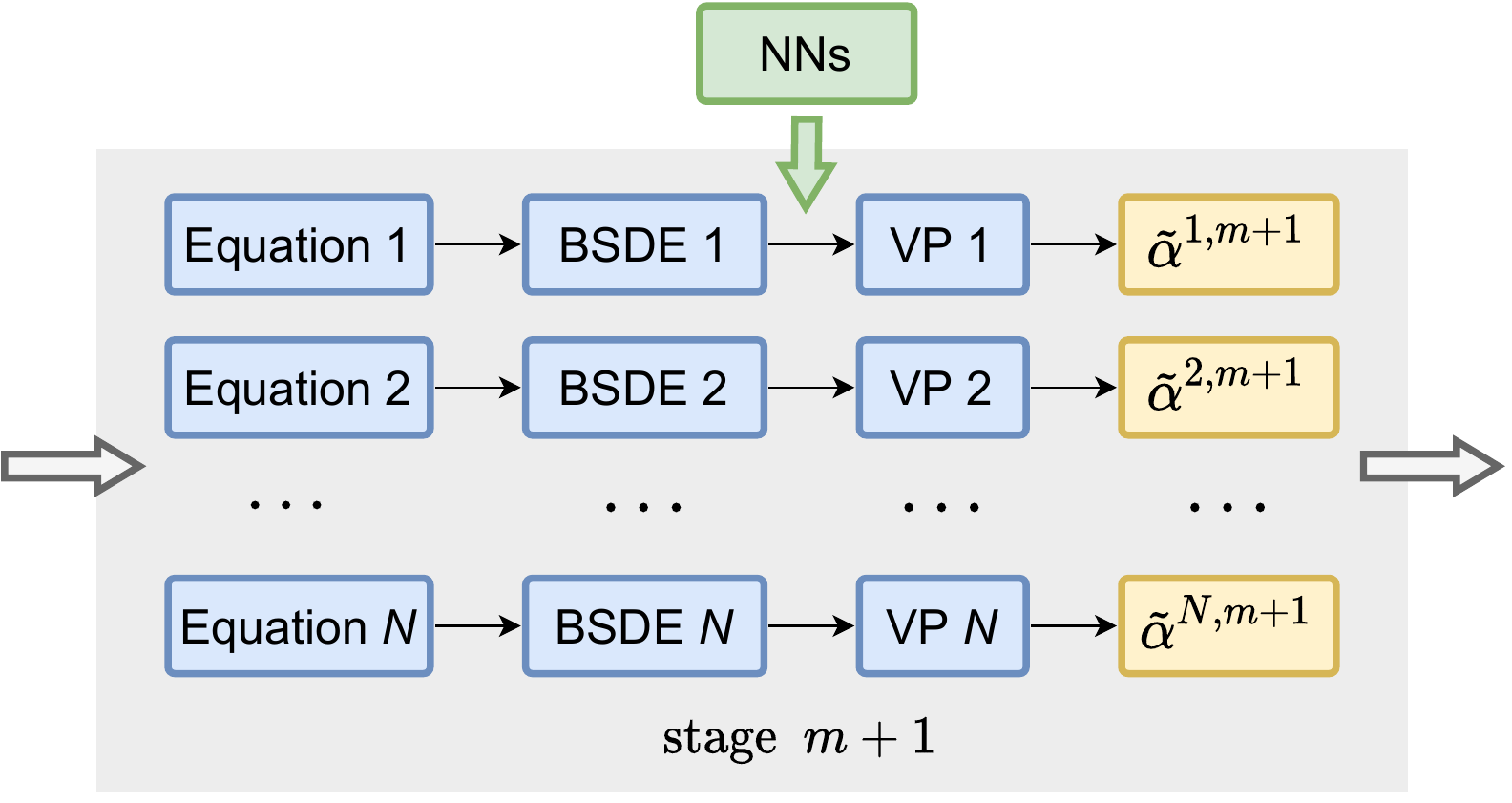}
  \caption{Illustration of one stage of \textit{enhanced deep fictitious play.} At the $(m+1)^{th}$ stage, one needs to solve the PDEs \eqref{def:HJB-DFP}, which is approximated by solving the BSDEs \eqref{eq:BSDE_forward}-\eqref{eq:BSDE_backward}. Then with the help of neural networks, one solves the variational problem (VP) given by Equation~\eqref{def:varitional_problem} to get the optimal strategy.}\label{fig:EDFP}
\end{figure}

We illustrate one stage of enhanced deep fictitious play in Figure~\ref{fig:EDFP}. At the $(m+1)^{th}$ stage, given the optimal policies $\bm{\alpha}^{m}$ at the previous stage, for $n = 1, \ldots, N$, the algorithm solves the following partial differential equations (PDEs),
\begin{multline}\label{def:HJB-DFP}
\partial_t V^{n, m+1} \\
+ \inf_{\alpha^n \in [0,1]^2} H^n(t,\bm x,(\alpha^n,\bm{\alpha}^{-n,m})(t, \bx), \nabla_{\bx} V^{n, m+1})\\ 
+ \half \text{Tr}(\Sigma(\bx)\transpose \text{Hess}_{\bx} V^{n, m+1} \Sigma(\bx)) = 0,
\end{multline} 
with $V^{n, m+1}(T,\bx) = 0,$
and obtains the optimal strategy of the $(m+1)^{th}$ stage:
\begin{align}\label{def:alphaast}
&\alpha^{n,m+1} = \\ &\argmin_{\alpha^n \in [0,1]^2} H^n(t, \bm x, (\alpha^n,\bm{\alpha}^{-n,m})(t, \bx), \nabla_{\bx} V^{n, m+1}(t, \bx)).
\end{align}



For simplicity of notations, we omit the stage number $m$ in the superscript in the following discussions. 
The solution to Equation~\eqref{def:HJB-DFP} is approximated by solving the equivalent backward stochastic differential equations (BSDEs) using neural networks \cite{MR3847747}:
\begin{empheq}[left=\empheqlbrace]{align}
    \bX_t^n =& \bx_0 + \int_{0}^{t} \mu^n(s, \bX_s^n; \bm{\alpha}^{-n}(s, \bX_s^n)) \ud s \\ &+ \int_{0}^{t}\Sigma(\bX_s^n)\, \mathrm{d}\bW_s, \label{eq:BSDE_forward} \\
    Y_t^n =&  \int_{t}^{T}g^n(s, \bX_s^n, Z_s^n; \bm{\alpha}^{-n}(s, \bX_s^n))\ud s \\ &- \int_{t}^{T}(Z_s^n)\transpose\, \mathrm{d}\bW_s. \label{eq:BSDE_backward}
\end{empheq}
The nonlinear Feynman-Kac formula \cite{MR1176785} yields:
\begin{equation}\label{eq:FK relation}
    Y_t^n = V^n(t, \bX_t^n) \quad\text{and}\quad Z_t^n = \Sigma(\bX_t^n)\transpose\nabla_{\bx} V^n(t, \bX_t^n).
\end{equation}
Here $\mu^n$ and $g^n$ are derived by rewriting \eqref{def:HJB-DFP} to
$\partial_t V^n + \half \text{Tr}(\Sigma(\bx)\transpose \text{Hess}_{\bx} V^n \Sigma(\bx)) + \mu^n(t, \bx; \balpha^-n)\cdot \nabla_{\bx} V^n+g^n(t, \bx, \Sigma(\bx)\transpose\nabla_{\bx} V^n; \balpha^{-n})=0$; see \cite[Appendix A.2]{xuan2020optimal}.
Notice that, we parametrized $V^n$ by neural networks (denote as $V$-networks) so $Y^n_t$ and $Z_t^n$ could all be computed by a function of $V$-networks. The loss function to update the $V$-network 
is constructed by simulating the BSDE along the time axis and penalizing the difference between the true terminal value and the simulated terminal value based on neural networks of $Y$. 


In Enhanced DFP, we further parameterize $\alpha^n$ (denote as $\alpha$-networks). In each stage, the loss function with respect to the $V$-network and the $\alpha$-network of player $n$ is defined by the weighted sum of two terms: the loss related to BSDE \eqref{eq:BSDE_forward}--\eqref{eq:BSDE_backward} to approximate its solution and the error of approximating the optimal strategy $\alpha^n$ by $\alpha$-networks. We denote $\|\cdot\|_2$ as the 2-norm,  $\alpha^n$ and  $\tilde{\alpha}^n$ as the derived and approximated optimal control of player $n$ in the current stage, $\tilde\balpha^{-n}=(\tilde\alpha^{1,m},...,\tilde {\alpha}^{n-1,m},\tilde \alpha^{n+1,m},...,\tilde{\alpha}^{N,m})$ as the collection of approximated optimal controls from the last stage except player $n$, and $\tau$ as a hyperparameter balancing the two types of errors in the loss function. Then the Enhanced DFP solves
 \begin{equation}
 \begin{split}\label{def:varitional_problem}
&\inf _{Y_{0}^{n},\tilde \alpha^n,\left\{Z_{t}^{n}\right\}_{0 \leq t \leq T}}\mathbb{E}(\left|Y_{T}^{n}\right|^{2}\\ &+\tau\int_0^T \ltwonorm{\alpha^n(s,\bX_s^n)-\tilde{\alpha}^n(s,\bX_s^n)}^2 \ud s) \\
\text { s.t. } \bX_t^n &= \bx_0 + \int_{0}^{t} \mu^n(s, \bX_s^n; \tilde\balpha^{-n}(s, \bX_s^n)) \ud s\\ & + \int_{0}^{t}\Sigma(\bX_s^n)\, \mathrm{d}\bW_s,  \\
    Y_t^n &=  Y_0^n-\int_{0}^{t}g^n(s, \bX_s^n, Z_s^n; \tilde\balpha^{-n}(s, \bX_s^n))\ud s\\ &+ \int_{0}^{t}(Z_s^n)\transpose\, \mathrm{d}\bW_s,\\
    {\alpha}^n(s,\bX_s^n)&=\argmin_{\beta^n} H^n(s, \bX_s^n, (\beta^n, \tilde\balpha^{-n})(s, \bX_s^n), Z_s^n).
\end{split}
\end{equation}
In each stage, there are two types of optimal strategies for player $n$: 1. the \textit{derived} optimal strategy $\alpha^{n}$ by solving $\argmin_{\beta^n} H^n$ in the last equation of \eqref{def:varitional_problem}; 2. the \textit{approximated} optimal strategy $\tilde{\alpha}^{n}$ also known as $\alpha$-networks for reducing the non-trivial cost of evaluating $\alpha^{n}$. Take stage $m+1$ as an example, $\alpha^{n,m+1}$ depends on players' last stage optimal policies $\bm{\alpha}^{-n, m}$ which in turn depends on $\bm{\alpha}^{-n, m-1}$. The evaluation of the current stage strategy $\alpha^{n,m+1}$ actually requires the recursive iteration of optimal strategies from all previous stages. Enhanced DFP unblocks the computation bottleneck by introducing approximated optimal strategy $\tilde{\alpha}^n$, which approximates $\alpha^n$. Although representing $\alpha^n$ with a neural network $\tilde{\alpha}^n$ introduces approximation errors, it allows us to efficiently access the proxy of the optimal strategy $\balpha^{-n}$ in the current stage by calling corresponding networks, instead of storing and calling all the previous strategies $\balpha^{-n, m-1}, \ldots, \balpha^{-n, 1}$ due to the recursive dependence. This is the key factor that {Enhanced Deep Fictitious Play} addresses leading to reduction in both time and memory complexity compared to {Deep Fictitious Play}.

To implement the loss function defined in \eqref{def:varitional_problem}, we discretize and simulate the BSDE by Euler's method with a partition on a time interval $[0,T]$. The expectation in the loss function is approximated by Monte Carlo samples of $\bW_s$ in the stochastic process. Then, we use the stochastic gradient descent method to update $V^n$ and $\tilde{\alpha}^{n}$ in the current stage for player $n$. In parallel, we update the $V$-network and $\alpha$-network for each player. The updated networks of each player will be observable for other players in future stages.

\section*{A Regional COVID-19 Study} 
In this section, we apply the multi-region stochastic SEIR model \eqref{def:St}--\eqref{def:J} to analyze optimal COVID-19 policies in three adjacent states: New York, New Jersey and Pennsylvania. This case study focuses on 180 days starting from 03/15/2020, and solves for the optimal policies of the three states corresponding to Nash Equilibrium by the Enhanced Deep Fictitious Play Algorithm. We denote New York (NY) as region 1, New Jersey (NJ) as region 2 and Pennsylvania (PA) as region 3, with population $P_1=19.54$ million, $P_2=8.91$ million, and $P_3=12.81$ million during the case study time range, respectively. We assume that (a) $90\%$  of any state's population resides in their own state at a given time; (b) the remaining population(travellers) visit the other states at an equal chance; (c) there is no travel outside of the three states, that is,  NY-NJ-PA is a closed system. The parameters in \eqref{def:St}-\eqref{def:J} are estimated based on the above assumptions and public information about COVID-19: $\beta = 0.17, \kappa = \frac{0.65\%}{13}, \lambda = \frac{1}{13}, \gamma = \frac{1}{5}, p = 228.7\times 10^{-5}, c = 73300/13$. Other parameters in the model are chosen at: $r = 0, w = 172.6, \chi = 1.96\times 10^6$. The hyperparameters, $\theta$ and $a$, which represent policy effectiveness and planners' views on the death of human beings will change the optimal policies. For results including vaccination controls, we point to \cite{OLIVARES2021106411}, which considers an optimal control problem for vaccines and testing of COVID-19. However, in the time period we study, vaccination was not available, so we ignore the health policy $h$  and mainly solve for the lockdown policy of each state. 

\begin{figure*}[ht]
    \centering
    \includegraphics[width = 0.75\textwidth, trim = {2em 2em 5em 6em}, clip, keepaspectratio=True]{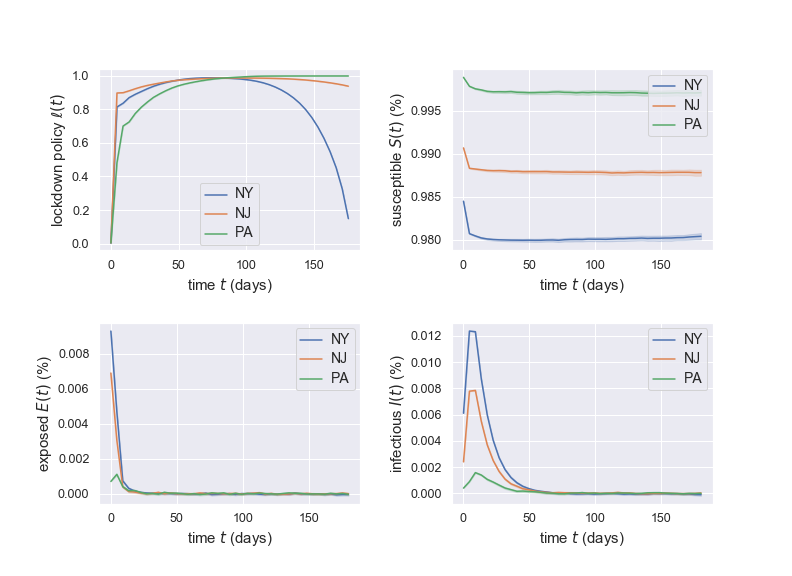}
    \caption{Plots of optimal policies (top-left), Susceptibles (top-right), Exposed (bottom-left) and Infectious (bottom-right) for three states: New York (blue), New Jersey (orange) and Pennsylvania (green). The shaded areas depict the mean and 95\% confidence interval over 256 sample paths. Choices of parameters are $a = 100$ and $\theta = 0.99$.}
    \label{fig:control}
\end{figure*}

Figure~\ref{fig:control} shows the Nash equilibrium policies in NY, NJ, and PA in a setting  where the policy effectiveness is $\theta=0.99$, {i.e.}, $99\%$ of the residents will follow the lockdown orders. The weight parameter quantifying each planner's view is $a=100$, {i.e.}, 
each governor values human life 100 times more than the economic value of a human life.
The resulting Nash equilibrium of this scenario corresponds to the planners taking action at an early stage by implementing strict lockdown policies and later relaxing the policy as the infections improve. In the end, the percentage of Susceptible, Exposed, Infectious, and Removed stays almost constant. The pandemic will be significantly mitigated in this scenario of proactive lockdown for both planners and residents. As a comparison, \cite[Figure 2]{xuan2020optimal} illustrates a scenario of how the pandemic gets out of control if governors show inaction or issue mild lockdown policies.




\section*{Acknowledgements}
The contents of this article are based on the authors' previous publication \cite{xuan2020optimal}. 
R.H.~was partially supported by the NSF grant DMS-1953035, the Faculty Career Development Award, the Research Assistance Program Award, the Early Career Faculty Acceleration funding, and the Regents' Junior Faculty Fellowship at University of California, Santa Barbara.
H.D.C. gratefully acknowledges partial support from the NSF grant DMS-1818821.

\bibliography{references}

\begin{bibdiv}
\begin{biblist}

\bib{alvarez2020simple}{techreport}{
      author={Alvarez, Fernando~E},
      author={Argente, David},
      author={Lippi, Francesco},
       title={A simple planning problem for covid-19 lockdown},
 institution={National Bureau of Economic Research},
        date={2020},
}

\bib{MR2428373}{incollection}{
      author={Allen, Linda J.~S.},
       title={An introduction to stochastic epidemic models},
        date={2008},
   booktitle={Mathematical epidemiology},
      series={Lecture Notes in Math.},
      volume={1945},
   publisher={Springer, Berlin},
       pages={81\ndash 130},
         url={https://doi.org/10.1007/978-3-540-78911-6_3},
      review={\MR{2428373}},
}

\bib{MR4121530}{article}{
      author={Abbasi, Zohreh},
      author={Zamani, Iman},
      author={Mehra, Amir Hossein~Amiri},
      author={Shafieirad, Mohsen},
      author={Ibeas, Asier},
       title={Optimal control design of impulsive {SQEIAR} epidemic models with
  application to {COVID}-19},
        date={2020},
        ISSN={0960-0779},
     journal={Chaos Solitons Fractals},
      volume={139},
       pages={110054, 20},
         url={https://doi.org/10.1016/j.chaos.2020.110054},
      review={\MR{4121530}},
}

\bib{MR3797719}{article}{
      author={Bottou, L\'{e}on},
      author={Curtis, Frank~E.},
      author={Nocedal, Jorge},
       title={Optimization methods for large-scale machine learning},
        date={2018},
        ISSN={0036-1445},
     journal={SIAM Rev.},
      volume={60},
      number={2},
       pages={223\ndash 311},
         url={https://doi.org/10.1137/16M1080173},
      review={\MR{3797719}},
}

\bib{bernoulli1760essai}{article}{
      author={Bernoulli, Daniel},
       title={Essai d'une nouvelle analyse de la mortalit{\'e} caus{\'e}e par
  la petite v{\'e}role, et des avantages de l'inoculation pour la
  pr{\'e}venir},
        date={1760},
     journal={Histoire de l'Acad., Roy. Sci.(Paris) avec Mem},
       pages={1\ndash 45},
}

\bib{Br:49}{techreport}{
      author={Brown, G.~W.},
       title={Some notes on computation of games solutions},
 institution={Rand Corp Santa Monica CA},
        date={1949},
}

\bib{MR0056265}{incollection}{
      author={Brown, George~W.},
       title={Iterative solution of games by fictitious play},
        date={1951},
   booktitle={Activity {A}nalysis of {P}roduction and {A}llocation},
      series={Cowles Commission Monograph No. 13},
   publisher={John Wiley \& Sons, Inc., New York, N.Y.; Chapman \& Hall, Ltd.,
  London},
       pages={374\ndash 376},
      review={\MR{0056265}},
}

\bib{MR3752669}{book}{
      author={Carmona, Ren\'{e}},
      author={Delarue, Fran\c{c}ois},
       title={Probabilistic theory of mean field games with applications. {I}},
      series={Probability Theory and Stochastic Modelling},
   publisher={Springer, Cham},
        date={2018},
      volume={83},
        ISBN={978-3-319-56437-1; 978-3-319-58920-6},
        note={Mean field FBSDEs, control, and games},
      review={\MR{3752669}},
}

\bib{cho2014learning}{article}{
      author={Cho, Kyunghyun},
      author={Van~Merri{\"e}nboer, Bart},
      author={Gulcehre, Caglar},
      author={Bahdanau, Dzmitry},
      author={Bougares, Fethi},
      author={Schwenk, Holger},
      author={Bengio, Yoshua},
       title={Learning phrase representations using rnn encoder-decoder for
  statistical machine translation},
        date={2014},
     journal={arXiv preprint arXiv:1406.1078},
}

\bib{grmek1997beginnings}{article}{
      author={Grmek, MD},
      author={Buchet, C},
       title={The beginnings of maritime quarantine},
        date={1997},
     journal={Man, health and the sea. Paris: Honor{\'e} Champion},
       pages={39\ndash 60},
}

\bib{han2016deep}{article}{
      author={Han, Jiequn},
      author={E, Weinan},
       title={Deep learning approximation for stochastic control problems},
        date={2016},
     journal={arXiv preprint arXiv:1611.07422},
}

\bib{HaHu:19}{inproceedings}{
      author={Han, J.},
      author={Hu, R.},
       title={Deep fictitious play for finding {Markovian} {Nash} equilibrium
  in multi-agent games},
        date={2020},
   booktitle={Proceedings of the first mathematical and scientific machine
  learning conference ({MSML})},
      volume={107},
       pages={107:221\ndash 245},
}

\bib{MR3847747}{article}{
      author={Han, Jiequn},
      author={Jentzen, Arnulf},
      author={E, Weinan},
       title={Solving high-dimensional partial differential equations using
  deep learning},
        date={2018},
        ISSN={0027-8424},
     journal={Proc. Natl. Acad. Sci. USA},
      volume={115},
      number={34},
       pages={8505\ndash 8510},
         url={https://doi.org/10.1073/pnas.1718942115},
      review={\MR{3847747}},
}

\bib{huLauriere2022}{article}{
      author={Hu, Ruimeng},
      author={Lauriere, Mathieu},
       title={Recent developments in machine learning methods for stochastic
  control and games},
        date={2022},
     journal={Available at SSRN: https://ssrn.com/abstract=4096569},
}

\bib{MR1025386}{article}{
      author={Kryscio, Richard~J.},
      author={Lef\`evre, Claude},
       title={On the extinction of the {S}-{I}-{S} stochastic logistic
  epidemic},
        date={1989},
        ISSN={0021-9002},
     journal={J. Appl. Probab.},
      volume={26},
      number={4},
       pages={685\ndash 694},
         url={https://doi.org/10.1017/s002190020002756x},
      review={\MR{1025386}},
}

\bib{lecun2015deep}{article}{
      author={LeCun, Yann},
      author={Bengio, Yoshua},
      author={Hinton, Geoffrey},
       title={Deep learning},
        date={2015},
     journal={nature},
      volume={521},
      number={7553},
       pages={436\ndash 444},
}

\bib{lecun1995convolutional}{article}{
      author={LeCun, Yann},
      author={Bengio, Yoshua},
      author={others},
       title={Convolutional networks for images, speech, and time series},
        date={1995},
     journal={The handbook of brain theory and neural networks},
      volume={3361},
      number={10},
       pages={1995},
}

\bib{MR908379}{article}{
      author={Liu, Wei~Min},
      author={Hethcote, Herbert~W.},
      author={Levin, Simon~A.},
       title={Dynamical behavior of epidemiological models with nonlinear
  incidence rates},
        date={1987},
        ISSN={0303-6812},
     journal={J. Math. Biol.},
      volume={25},
      number={4},
       pages={359\ndash 380},
         url={https://doi.org/10.1007/BF00277162},
      review={\MR{908379}},
}

\bib{MR43432}{article}{
      author={Nash, John},
       title={Non-cooperative games},
        date={1951},
        ISSN={0003-486X},
     journal={Ann. of Math. (2)},
      volume={54},
       pages={286\ndash 295},
         url={https://doi.org/10.2307/1969529},
      review={\MR{43432}},
}

\bib{nordstrom2022risk}{article}{
      author={Nordstr{\"o}m, Peter},
      author={Ballin, Marcel},
      author={Nordstr{\"o}m, Anna},
       title={Risk of sars-cov-2 reinfection and covid-19 hospitalisation in
  individuals with natural and hybrid immunity: a retrospective, total
  population cohort study in sweden},
        date={2022},
     journal={The Lancet Infectious Diseases},
      volume={22},
      number={6},
       pages={781\ndash 790},
}

\bib{OLIVARES2021106411}{article}{
      author={Olivares, Alberto},
      author={Staffetti, Ernesto},
       title={Optimal control-based vaccination and testing strategies for
  covid-19},
        date={2021},
        ISSN={0169-2607},
     journal={Computer Methods and Programs in Biomedicine},
      volume={211},
       pages={106411},
  url={https://www.sciencedirect.com/science/article/pii/S0169260721004855},
}

\bib{MR1176785}{incollection}{
      author={Pardoux, \'{E}.},
      author={Peng, S.},
       title={Backward stochastic differential equations and quasilinear
  parabolic partial differential equations},
        date={1992},
   booktitle={Stochastic partial differential equations and their applications
  ({C}harlotte, {NC}, 1991)},
      series={Lect. Notes Control Inf. Sci.},
      volume={176},
   publisher={Springer, Berlin},
       pages={200\ndash 217},
         url={https://doi.org/10.1007/BFb0007334},
      review={\MR{1176785}},
}

\bib{tornatore2005stability}{article}{
      author={Tornatore, Elisabetta},
      author={Buccellato, Stefania~Maria},
      author={Vetro, Pasquale},
       title={Stability of a stochastic {SIR} system},
        date={2005},
     journal={Physica A: Statistical Mechanics and its Applications},
      volume={354},
       pages={111\ndash 126},
}

\bib{vaswani2017attention}{article}{
      author={Vaswani, Ashish},
      author={Shazeer, Noam},
      author={Parmar, Niki},
      author={Uszkoreit, Jakob},
      author={Jones, Llion},
      author={Gomez, Aidan~N},
      author={Kaiser, {\L}ukasz},
      author={Polosukhin, Illia},
       title={Attention is all you need},
        date={2017},
     journal={Advances in neural information processing systems},
      volume={30},
}

\bib{xuan2020optimal}{inproceedings}{
      author={Xuan, Yao},
      author={Balkin, Robert},
      author={Han, Jiequn},
      author={Hu, Ruimeng},
      author={Ceniceros, Hector~D},
       title={Optimal policies for a pandemic: a stochastic game approach and a
  deep learning algorithm},
organization={PMLR},
        date={2022},
   booktitle={Proceedings of the second mathematical and scientific machine
  learning conference ({MSML})},
       pages={145:987\ndash 1012},
}

\end{biblist}
\end{bibdiv}

\end{document}